\documentclass[mathpazo]{cicp}

\usepackage{amsmath, amsthm,amsfonts}
\usepackage{tikz}
\usepackage{xcolor}
\usepackage{float}
\restylefloat{figure}
\newcommand{\Real}{\mathbb{R}}
\newcommand{\A}{\mathcal{A}}

\newcommand{\C}{\mathcal{C}}

\DeclareMathOperator{\e}{e}

\newcommand{\Pl}{\mathbb{P}}
\newcommand{\B}{\mathcal{B}}

\begin{document}

\title{Adapted nested force-gradient integrators: the Schwinger model case  }

\author[Shcherbakov D  et.~al.]{Dmitry Shcherbakov\affil{1}\comma\corrauth,
      Matthias Ehrhardt\affil{1},  Jacob Finkenrath\affil{2}, Michael G\"unther\affil{1},
      Francesco Knechtli\affil{3} and Michael Peardon\affil{4}}
\address{\affilnum{1}\ Lehrstuhl f\"ur Angewandte Mathematik und Numerische Analysis,
Bergische Universit\"{a}t Wuppertal,
Gau{\ss}strasse 20, 42119 Wuppertal, Germany \\
          \affilnum{2}\  CaSToRC, CyI,
        20 Constantinou Kavafi Street, 2121 Nicosia, Cyprus  \\
          \affilnum{3}\  Theoretische Physik,
Bergische Universit\"{a}t Wuppertal,
Gau{\ss}strasse 20, 42119 Wuppertal, Germany \\
          \affilnum{4}\  School of Mathematics, Trinity College, Dublin 2, Ireland }
\emails{{\tt shcherbakov@math.uni-wuppertal.de} (D.~Shcherbakov), {\tt ehrhardt@math.uni-wuppertal.de} (M.~Ehrhardt),{\tt j.finkenrath@cyi.ac.cy} (J.~Finkenrath),
         {\tt guenther@math.uni-wuppertal.de} (M.~G\"unther), 
         {\tt knechtli@physik.uni-wuppertal.de} (F.~Knechtli), {\tt mjp@maths.tcd.ie} (M.~Peardon)}

% 
% \author{Dmitry Shcherbakov\corref{cor1}} 
% \ead{shcherbakov@math.uni-wuppertal.de}
% \cortext[cor1]{corresponding author}
% 
% \author{ Matthias Ehrhardt}
% \ead{ehrhardt@math.uni-wuppertal.de}
% 
% \author{Michael G\"unther}
% \ead{guenther@math.uni-wuppertal.de}
% 
% \author{Jacob Finkenrath}
% \ead{finkenrath@physik.uni-wuppertal.de}
% 
% \author{Francesco Knechtli}
% \ead{knechtli@physik.uni-wuppertal.de}
% 
% \address{Lehrstuhl f\"ur Angewandte Mathematik und Numerische Analysis, Bergische Universit\"{a}t Wuppertal,
% Gau{\ss}strasse 20, 42119 Wuppertal, Germany}
% 
% \author{Michael Peardon}
% \ead{mjp@maths.tcd.ie}
% \address{School of Mathematics, Trinity College, Dublin 2, Ireland}

%%%%%%%%%%%%%%%%%%%%%%%%%%%%%%%%%%%%%%%%%%%%%%%%%%%
\begin{abstract}
We study a novel class of numerical integrators, the adapted nested 
force-gradient schemes, used within the molecular dynamics step of
the Hybrid Monte Carlo (HMC) algorithm.  We test these methods  
in the Schwinger model on the lattice, a well known benchmark problem.
We derive the analytical basis of nested force-gradient type methods 
and demonstrate the advantage of the proposed approach, namely reduced 
computational costs compared with other numerical integration schemes in HMC. 
\end{abstract}

\maketitle

%%%%%%%%%%%%%%%%%%%%%%%%%%%%%%%%%%%%%%%%%%%%%%%%%%%%%%%%%%%%%%%%%%%%%
\section{Introduction}

For the Hybrid Monte Carlo algorithm (HMC)\cite{HMC}, often used to study 
quantum chromodynamics (QCD) on the lattice, one is interested in efficient 
numerical time integration schemes which are optimal in terms of computational 
costs per trajectory for a given acceptance rate. High order numerical methods 
allow the use of larger step sizes, but demand a larger computational effort 
per step; low order schemes do not require such large computational costs per 
step, but need more steps per trajectory. So there is a need to balance these
opposing effects. 

Omelyan integration schemes \cite{OmMrFo03} of a force-gradient type have 
proved to be an efficient choice, since it is easy to obtain higher order
schemes that demand a small additional computational effort. These 
schemes use higher-order information from force-gradient
terms to both increase the convergence of the method and decrease the size of 
the leading error coefficient. 
Other ideas to achieve better efficiency for numerical time integrators are 
given by multirate or nested approaches. These schemes do not increase the 
order but reduce the computational costs per path  by recognizing the 
different dynamical time-scales generated by different parts of the action. 
Slow forces, 
which are usually expensive to evaluate, need only to be sampled at 
low frequency while fast forces which are usually cheap to evaluate need to be 
sampled at a high frequency.
A natural way to inherit the  advantages from both force-gradient type schemes 
and multirate approaches would be to combine these two ideas.

Previously, we studied the behavior of the adapted nested force-gradient scheme 
for the example of the $n$-body problem \cite{ShEhGuP15} and would like to learn
more about their usefulness for lattice field theory calculations. Due to the 
huge computational effort required for QCD simulations, it is 
natural to attempt an intermediate step first.
We chose the model problem of quantum electrodynamics (QED) in two dimensions, 
the Schwinger model \cite{Sch67}, since it is well-suited as a test 
case for new concepts and ideas which can be subsequently applied to more 
computationally demanding 
problems \cite{ChJaNaPo05}. As a lattice quantum field theory, it has 
many of the properties of more sophisticated models such as QCD, for example 
the numerical cost is still dominated by the fermion part of the action. 
The fact that this model, with far fewer degrees of freedom, does not require 
such large computational effort makes it the perfect choice for testing purposes.

We compare the behavior of numerical time integration schemes currently used
for HMC \cite{OmMrFo03} with the nested force-gradient integrator 
\cite{KeSipCl10} and the adapted version introduced in  \cite{ShEhGuP15}. We 
investigate the computational costs needed to perform numerical calculations, 
as well as the effort required to achieve a satisfactory 
acceptance rate during the HMC evolution. Our goal is to find a numerical 
scheme for the HMC algorithm which would provide a sufficiently high 
acceptance rate while not drastically increasing the simulation time. 

 The paper is organized as follows. In Section 2 we give a short overview of 
the HMC algorithm  and numerical schemes for time integration, which are used 
in HMC. In Section 3 we present the  2-dimensional  Schwinger model and 
introduce the idea of the force-gradient approach  and the resulting novel  
class of numerical schemes.  Section 4 is devoted to the results of a 
comparison between widely used algorithms and the new approach and Section 5 
draws our conclusion. 

% % % % % % % % % % % % % % % % % % % % % % % % % % % % % % % % % % % % % % 

\section{Geometric integrators for HMC}
In this section we provide a general overview of the HMC algorithm \cite{HMC}
to introduce our novel integrator. We also present some standard numerical time 
integrating methods used in HMC, as well state-of-the-art numerical schemes, 
which we later compare by applying them to the two-dimensional Schwinger model. 

\subsection{Overview of HMC}

In the Hybrid Monte Carlo algorithm, the quantum lattice field theory is 
embedded in a higher-dimensional classical system through the introduction of a fictitious (simulation) time \cite{HMC}. 
The gauge field $U$ is associated 
with its (fictitious) conjugate momenta $P$, and the classical system is 
described by the Hamiltonian,
\begin{equation} \label{eq:ham}
 H = \A[P] + \B[U],
\end{equation}
where $\A[P]$ and $\B[U]$ represent the kinetic and potential energy 
respectively. 

For a given configuration $U$, a new configuration $U'$ is generated by 
performing an HMC update $ U \rightarrow U'$, which consists of two steps:
\begin{itemize}
 \item \textbf{Molecular Dynamics trajectory:}
 Evolve the gauge fields $U$, elements of a Lie group,
and the momenta $P$, elements of the corresponding Lie algebra,   
in a fictitious computer time $t$ according to Hamilton's equations of motions
 \begin{align}\label{eq:eq_motion}
 \dot{P} = -\frac{\partial H}{\partial U} = -F_{V}(U), &&  \dot{U} = P U.  
\end{align}
Since analytical solutions are not available in general, numerical methods 
must be used to solve the system of Eqn.~\eqref{eq:eq_motion}. 
The discrete updates of $U$ and $P$ with an integration step $h$ are 
\begin{align*}
& \e^{\A h}: ~~ U ~~ \rightarrow ~~ U' = \exp(i P h )U  \\
& \e^{\B h}: ~~ P ~~ \rightarrow ~~ P' = P - i h  F_{V}(U),
\end{align*}
leading to a first-order approximation at time $t+h$.
Since the momenta $P$ are elements of Lie algebra, we have an additive update 
of $P$. On the other hand, the links $U$ must be  elements of  the Lie group, 
therefore an exponential update is used for $U$ to preserve the underlying 
group structure.
 \item \textbf{Metropolis step:}
 Accept or reject the new configuration $(U',P')$ with probability
 $$\mathcal{P}(U \rightarrow U') = \min\left( 1, \e^{-\Delta H}\right), $$
 where $\Delta H = H(U',P')-H(U,P)$.
\end{itemize}

\subsection{Integrators used for Molecular Dynamics}
In this paper we are concerned with numerical time integration schemes, which 
preserve the fundamental properties of geometric integration, 
time-reversibility and volume-pre\-servation. 
All numerical schemes presented below possess these necessary properties.

{\bf Basic schemes:}
Well-known, commonly used integration schemes in molecular dynamics are given by 
\begin{itemize}
 \item the leap-frog method, a 3-stage composition scheme of the discrete updates defined above:
 \begin{equation}\label{eq:leap}
\Delta(h)=\e^{h\frac{\hat{\B}}{2}}\e^{h\hat{\A}}\e^{h\frac{\hat{\B}}{2}},
\end{equation}
\item and a 5-stage extension widely used in QCD computations:
\begin{equation}\label{eq:5stage}
\Delta_5(h) = \e^{ \frac{1}{6} h \hat{\B} } \e^{ \frac{1}{2} h \hat{\A} }  \e ^{ \frac{2}{3} h \hat{\B} } \e^{ \frac{1}{2} h \hat{\A} } \e^{ \frac{1}{6} h \hat{\B}}.
\end{equation}
\end{itemize}

{\bf Force gradient schemes:}
Force-gradient schemes increase accuracy by using additional information from 
the force gradient term $\C= \{\B, \{  \A ,\B \}\}$, with $\{~,~\}$ defining Lie brackets. The  5-stage force-gradient scheme proposed by Omelyan et al 
\cite{OmMrFo03} is the simplest;
\begin{equation}\label{eq:5stage_fg}
 \Delta_{5C}(h) = \e^{ \frac{1}{6} h \hat{\B} } \e^{ \frac{1}{2} h \hat{\A} }  \e ^{ \frac{2}{3} h \hat{\B} - \frac{1}{72} h^{3} \C } 
  \e^{\frac{1}{2} h \hat{\A} }  \e^{ \frac{1}{6} h \hat{\B} }.
  \end{equation}
   Here we also test the modification of the force-gradient 
method \eqref{eq:5stage_fg} proposed in \cite{YiMa}, where the force-gradient term $\C$ is approximated via a Taylor expansion.
An extension is given by the 11-stage decomposition \cite{OmMrFo03}, recently 
implemented as the integrator in the open source code 
openQCD  as one of  the standard options \cite{LS13} \footnotesize
 \begin{equation}\label{eq:11_stage}
   \Delta_{11}(h) = \e^{ \sigma h \hat{\B} } \e^{\eta h \hat{\A} }  \e ^{ \lambda h \hat{\B}  } 
  \e^{\theta h \hat{\A} }  \e^{ (1-2(\lambda+\sigma)) \frac{h}{2} \hat{\B} }\e^{(1-2(\theta+\eta)) h \hat{\A} }  \e ^{ (1-2(\lambda+\sigma)) \frac{h}{2} \hat{\B} } 
  \e^{\theta h \hat{\A} }  \e^{ \lambda h \hat{\B} }\e^{\eta h \hat{\A} }  \e^{ \sigma h \hat{\B} },
  \end{equation}\normalsize
where  $\sigma$, $\theta$,  $\lambda$ and $\eta$ are parameters from equation 
(71) in Ref.~\cite{OmMrFo03}.

{\bf Nested Schemes:} QED and QCD problems usually lead to Hamiltonians with 
the following fine structure
\begin{equation} \label{eq:multiscale}
 H = \A[P] + \B_{1}[U]+ \B_{2}[U],
\end{equation}
where the action of the system can be split into two parts: a  fast action 
$\B_{1}$ such as the gauge action, and a slow part $\B_{2}$, for example, the  
fermion action. This allows us to apply the idea of  multirate  schemes (an 
idea known as nested integration in physics literature)\cite{Sexton92} in order to 
reduce the computational effort.
At first we consider the nested version of the leap-frog method \eqref{eq:leap} 
\begin{equation}\label{eq:leap_multi}
\hat{\Delta}(h) = 
    \e^{  \frac{h}{2} \hat{\B}_{2}} \Delta_1\left( h \right)_{M}  \e^{  \frac{h}{2} \hat{\B}_{2}},
\end{equation}
where the inner cheaper system $H=\A[P]+\B_{1}[U]$ is solved by 
\begin{equation*}
\Delta_1(h)_{M} =  \left(\e^ {\frac{h}{2M}\hat{\B}_{1}}  \e^{\frac{h}{M}\hat{\A}}  \e^{\frac{h}{2M}\hat{\B}_{1}}\right)^{M}, 
\end{equation*}
with $M$ being a number of iterations for the fast part of the action. Our main goal is to compare 
the above-mentioned methods with more elaborated nested schemes: in \cite{ShEhGuP15}, a similar 5-stage decomposition scheme has been recently introduced: 
\begin{equation}\label{eq:5stage_multi}
\hat{\Delta}_5(h) = 
 \e^{ \frac{1}{6} h \hat{\B}_{2}}  \Delta_1\left( \frac{h}{2} \right)_{M}  \e^{  \frac{2}{3} h \hat{\B}_{2}}   
\Delta_1\left(\frac{h}{2}\right)_{M} \e^{\frac{1}{6} h\hat{\B}_{2} }.
\end{equation}

A nested version of \eqref{eq:5stage_fg}, which has been used in~\cite{Bazavovetal2012} reads 
%is proposed in \cite{KeSipCl10}
\begin{equation}\label{eq:5stage_multi_fg}
\hat{\Delta}_{5C}(h) =   \e ^{\frac{1}{6} h \hat{\B}_{2}} \Delta_2\left( \frac{h}{2} \right)_{M}
\e^{ \frac{2}{3} h \hat{\B}_{2} + 
\frac{1}{72}h^{3}\C_f}    \Delta_2\left(\frac{h}{2}\right)_{M}  \e^{\frac{1}{6} h\hat{\B}_{2}}, 
\end{equation}
where 
\begin{equation*}
\Delta_2(h)_{M} = \left(\e^{ \frac{1}{6M} h \hat{\B_1} } \e^{ \frac{1}{2M} h \hat{\A} } 
\e ^{ \frac{2}{3M} h \hat{\B_1} + \frac{1}{72} \left(\frac{h}{M}\right)^{3} \C_g } 
  \e^{\frac{1}{2M} h \hat{\A} }  \e^{ \frac{1}{6M} h \hat{\B_1} }\right)^{M}, 
  \end{equation*}
with $\C_g = \{\B_1,\{\A,\B_1\}\}$ and $\C_f= \{\B_2,\{\A,\B_2\}\}$. 
In the limit $M\to\infty$ we have $\Delta_1=\Delta_2$.
Note that 
this approach uses force-gradient information at 
all levels, i.e., the high computational cost of high order schemes appears 
at all levels.

One may overcome this problem by 
using schemes of different order at the different levels without 
losing the effective high order of the overall multirate scheme. For the 
latter, we include appropriate force gradient information
as we explain in the following 
%MG 31.8
%.
for the case of 
%For this, we restricted the system to 
two time levels, 
%as the 
%ultimate goal will be to develop nested force gradient schemes for %lattice
%QCD, 
where the gauge action plays the role of the fast and cheap part, 
and the fermionic action plays the role of the slow and expensive part.
%, i.e., only two time scales are involved. 
%MG 31.8
The reasoning is as follows: if one uses the 5-stage Sexton-Weingarten integrator of second order for the slow action, 
and approximates the fast action by $m$ Leap-frog steps of step size $h/(2m)$, 
%If one uses a second-order scheme for 
%the fast action with time step $h/m$ and a fourth-order scheme for the %slow 
%action with time step $h$, 
the error of the overall multirate scheme will be 
of order ${\cal O}(h^2)+{\cal O}((\frac{h}{m})^2)+{\cal O}(h^4)$. With the use 
of force gradient information only at the slowest level it is possible to 
cancel the leading error term of order ${\cal O}(h^2)$. As $m \geq \frac{1}{h}$
usually holds in the multirate setting, the overall order is then given by the 
leading error term of order ${\cal O}(h^4)$, i.e., the scheme has an effective 
order of four. 
One example for such a scheme for problems of type \eqref{eq:multiscale} is given by the  5-stage nested force-gradient scheme
introduced in \cite{ShEhGuP15}
\begin{equation}\label{eq:ad_5stage_multi_fg}
\tilde{\Delta}_{5C}(h) =   \e ^{\frac{1}{6} h \hat{\B}_{2}} \Delta_1\left( \frac{h}{2} \right)_{M} 
\e^{ \frac{2}{3} h \hat{\B}_{2} + 
\frac{1}{72}h^{3}\C_f}    \Delta_1\left(\frac{h}{2}\right)_{M}  \e^{\frac{1}{6} h\hat{\B}_{2}}. 
\end{equation}
To summarize, the adapted scheme~\eqref{eq:ad_5stage_multi_fg} differs from 
the original one ~\eqref{eq:5stage_multi_fg} in two 
perspectives:
\begin{itemize}
\item The force gradient scheme for the fast action is replaced by a leap-frog 
scheme.
\item Only the part $\{\B_2,\{\A,\B_2\}\}$ of the full force gradient
$ \{\B_1+\B_2,\{\A,\B_1+\B_2\}\}$ is needed to gain the effective order of four.
\end{itemize}

The numerical 
schemes \eqref{eq:leap}-\eqref{eq:5stage} and \eqref{eq:leap_multi}-\eqref{eq:5stage_multi} are  second order 
 convergent schemes. Methods \eqref{eq:5stage_fg}-\eqref{eq:11_stage} and  \eqref{eq:5stage_multi_fg} - \eqref{eq:ad_5stage_multi_fg}
 have the fourth  order of  convergence. We do not consider integrators of  higher order than four since the computational costs 
 are too high. The schemes of the  same convergence order differ from each other by the number of stages
(updates of momenta and links per time step). Usually methods with more stages have  smaller leading error coefficients and therefore 
have better accuracy, but higher computational costs.  We would like to 
determine which  integrator would represent the best compromise
between high accuracy and computational efficiency.

We will apply all these numerical integration schemes \eqref{eq:leap}--\eqref{eq:ad_5stage_multi_fg}  
to the two-dimensional Schwinger model. The most challenging task 
from the theoretical point of view is to derive the force-gradient term $\C$. In the next section 
we introduce the Schwinger model and explain how to obtain the force-gradient term.

\section{The Schwinger model and its force gradient terms}

The 2 dimensional Schwinger model is defined by the following Hamiltonian function
\begin{equation}\label{eq:schwinger}
{ H = \frac{1}{2}\sum_{n=1,\mu=1}^{V,2} p_{n,\mu}^2 + S_{full}[U] = \frac{1}{2}\sum_{n=1,\mu=1}^{V,2} p_{n,\mu}^2 + S_G[U] +S_F[U].}
\end{equation}
{with $V=L\times T$ the volume of the lattice.}
Unlike QCD, where $U \in SU(3)$ and $p_{n,\mu} \in su(3)$, for this QED problem \eqref{eq:schwinger}, 
the links   $U$ are the elements of the Lie group $U(1)$ and the momenta $p_{n,\mu}$ 
belong to $\Real$, which represents the Lie algebra of the  group $U(1)$. This 
makes this test example \eqref{eq:schwinger} very cheap in  terms of the computational time. This together with the fact that the Schwinger model also shares
many of the features of QCD simulations, makes the Schwinger model an excellent
test example when considering numerical integrators:
a fast dynamics given by the computationally cheap
gauge part $S_G[U]$ of the action demanding small step sizes, and a slow dynamics given by the 
computationally expensive 
 fermion part $S_F[U]$ allowing large step sizes. 

The pure gauge part of the action $S_G$ {sums up over all 
plaquettes $\Pl(n)$ in the two-dimensional lattice with}
$$\Pl(n) =  U_{1}(n) U_{2}(n+\hat{1}) U^{\dagger}_{1}(n+\hat{2})U^{\dagger}_{2}(n),$$ 
{and is given by}
\begin{equation}\label{eq:gauge_action}
{ S_G = \beta \sum\limits_{n=1}^V \left( 1 -  \operatorname{Re}  \Pl(n)  \right).}
\end{equation}
{The links $U$ can be written in the form  $U_{\mu}(n) = e^{i q_\mu(n)} \in \textrm{U}(1)$} and connect the 
sites $n$ and $n+\hat{\mu}$ on the lattice;
 $q_\mu(n)\in \left[ - \pi, \pi  \right] $, $\mu$, $\nu$ $\in \{x,t\}$ are respectively space and time directions and 
 $\beta$ is a coupling constant. {Note that from now on we will set the lattice spacing $a=1$.}
 
 The fermion part of the action $S_F$ is  given by  
 \begin{equation}\label{eq:fer_action}
  S_F = \eta^{\dagger} \left(  D^{\dagger}D\right)^{-1} \eta,
 \end{equation}
{ where $\eta$ is a complex pseudofermion field.
 Here,  $D$ denotes the Wilson--Dirac operator given by}
\begin{equation*}
{ D_{n,m} = (2+m_0)\delta_{n,m} 
  - \frac{1}{2} \sum\limits_{\mu=1}^2 \left(  
(1-\sigma_{\mu}) U_{\mu}(n) \delta_{n,m-\hat{\mu}} + 
(1+\sigma_{\mu}) U^{\dagger}_{\mu}(n-\hat{\mu}) \delta_{n,m+\hat{\mu}}  \right), }
\end{equation*}
where   $\sigma_{\mu}$ are the Pauli matrices 
\begin{equation*}
 \sigma_{1} = \begin{pmatrix}
               0 & 1\\
               1 & 0
              \end{pmatrix}~~ \text{and} ~~            
 \sigma_{2} = \begin{pmatrix}
               0 & -i\\
               i & 0
              \end{pmatrix}.         
\end{equation*}
${ m_0}$ is {the mass parameter and the Kronecker delta} $\delta_{n,m}$ {acts on
the pseudofermion field by} $\sum_{m=1}^V \delta_{n,m} \eta(m) = \eta(n)$ {with $\eta(n)$ the pseudofermion field, 
a vector in the two-dimensional spinor space taking values at each lattice 
point $n$.}
In order to proceed with the numerical integration we need to obtain the force 
$F$ and the force gradient term $\C$. 
The force term {$F(n,\mu)$ with respect to the link $U_\mu(n)$} is given by the first derivative of the action $S_{full}$ and can be written as 
\begin{equation}\label{eq:force}
{F(n,\mu)= F_{S_G}(n,\mu) + F_{S_F}(n,\mu) = \frac{\partial S_G}{\partial q_{\mu}(n)} + \frac{\partial S_F}{\partial q_\mu(n)}.}
 \end{equation}
Since the numerical schemes \eqref{eq:5stage_multi}--\eqref{eq:ad_5stage_multi_fg} use the multi-rate approach,
the shifts in the momenta updates are split on $F_{S_G}$  and  $F_{S_F}$ and  we can consider them separately. The force terms
$F_{S_G}$  and  $F_{S_F}$ are obtained by 
differentiation over $U(1)$ group elements, which for the Schwinger model is 
the standard differentiation. 

The force associated with link $U_{\mu}(n)$ from the gauge action is given by
\begin{equation}\label{eq:force_gauge}
 { \beta g(n,\mu) := F_{S_G}(n,\mu)= \beta \left. \operatorname{Im}\left( \Pl(n) - \Pl(n-\hat{\nu}) \right) \right|_{\mu\neq \nu} \,.}
\end{equation}
The force term of the fermion part is given by
\small
  \begin{equation}\label{eq:force_fermion}
f(n,\mu) := F_{S_F} = -\operatorname{Im} \left[
\chi^{\dagger}(n) (1-\sigma_{\mu}) U_{\mu}(n) \xi(n+\hat{\mu}) -
\chi^{\dagger}(n+\hat{\mu}) (1+\sigma_{\mu}) U^{\dagger}_{\mu}(n) \xi(n)\right] \,
    \end{equation}
\normalsize
    where vectors $ \chi$ and $\xi$ are given
  \begin{align}\label{eq:chi_xi_vectors}
  \chi = D^{\dagger^{-1}} \eta, && \xi= D^{-1}D^{\dagger^{-1}} \eta.
  \end{align}

 For the numerical methods \eqref{eq:5stage_fg} and \eqref{eq:5stage_multi_fg} we need to find the force gradient term {$\C(n,\mu)$ with respect to the link $U_\mu(n)$}. 
 In case of the Schwinger model \eqref{eq:ham} this term reads
 \begin{equation}\label{eq:force_gradient}
  {\C(n,\mu) =  
  2 \sum\limits_{m=1,\nu=1}^{V,2} \frac{\partial S_{full}}{\partial q_{\nu}(m)} \frac{\partial^{2} S_{full}}{\partial q_{\nu}(m) \partial q_{\mu}(n)}.}
 \end{equation}

For simplicity we decompose the force gradient term \eqref{eq:force_gradient} in four parts 
 \begin{equation}\label{eq:force_grad_de}
 \begin{aligned}
  &{\C_{GG}=   2 \sum\limits_{m=1,\nu=1}^{V,2} \frac{\partial S_{G}}{\partial q_{\nu}(m)} \frac{\partial^{2} S_{G}}{\partial q_{\nu}(m) \partial q_{\mu}(n)},} 
  &{\C_{FG}=   2 \sum\limits_{m=1,\nu=1}^{V,2} \frac{\partial S_{F}}{\partial q_{\nu}(m)} \frac{\partial^{2} S_{G}}{\partial q_{\nu}(m) \partial q_{\mu}(n)},} \\
  &{\C_{GF}=   2\sum\limits_{m=1,\nu=1}^{V,2} \frac{\partial S_{G}}{\partial q_{\nu}(m)} \frac{\partial^{2} S_{F}}{\partial q_{\nu}(m) \partial q_{\mu}(n)}, }
  &{\C_{FF}=   2 \sum\limits_{m=1,\nu=1}^{V,2} \frac{\partial S_{F}}{\partial q_{\nu}(m)} \frac{\partial^{2} S_{F}}{\partial q_{\nu}(m) \partial q_{\mu}(n)}.}
  \end{aligned}
  \end{equation}
This decomposition is also useful since  the numerical integrator 
\eqref{eq:5stage_multi_fg} only uses the term ${\C_{FF}}$
by construction. As shown in \cite{ShEhGuP15}, to obtain the fourth order convergent
scheme \eqref{eq:5stage_multi_fg} from the second order  convergent method \eqref{eq:5stage_multi}
we must eliminate the leading error term, which is exactly represented by ${\C_{FF}}$.
For completeness we discuss all 4 parts below.

The ${\C_{GG}}$ part of the force-gradient term  is
 \normalsize
\begin{equation*}
  \begin{aligned}
\C_{GG} &= 2 \beta^{2}
 \left[ 
 \operatorname{Im}(4 \Pl_{1}(n,\mu)- \Pl_{2}(n,\mu) -\Pl_{3}(n,\mu) -\Pl_{4}(n,\mu) -\Pl_{5}(n,\mu)) \operatorname{Re}(\Pl_{1}(n,\mu)) \right. \\
& \left. -   \operatorname{Im}(4 \Pl_{2}(n,\mu) - \Pl_{1}(n,\mu) - \Pl_{6}(n,\mu) - \Pl_{7}(n,\mu) -\Pl_{8}(n,\mu))      \operatorname{Re}(\Pl_{2}(n,\mu))  \right]
\end{aligned}
\end{equation*}
\normalsize
with the set of plaquettes 
\begin{align*}
\Pl_{1}(n,\mu)
      & = U_{\mu}(n)            U_{\nu}(n+\hat{\mu}) U^{\dagger}_{\mu}(n+\hat{\nu})U^{\dagger}_{\nu}(n),\\
\Pl_{2}(n,\mu)
      & =  U_{\mu}(n-\hat{\nu}) U_{\nu}(n-\hat{\nu}+\hat{\mu}) U^{\dagger}_{\mu}(n)U^{\dagger}_{\nu}(n-\hat{\nu}),\\
\Pl_{3} (n,\mu)
      & = U_{\mu}(n+\hat{\mu}) U_{\nu}(n+2\hat{\mu}) U^{\dagger}_{\mu}(n+\hat{\nu}+\hat{\mu})U^{\dagger}_{\nu}(n+\hat{\mu}),\\ 
\Pl_{4}(n,\mu)
      & = U_{\mu}(n+\hat{\nu}) U_{\nu}(n+\hat{\mu}+\hat{\nu}) U^{\dagger}_{\mu}(n+2\hat{\nu})U^{\dagger}_{\nu}(n+\hat{\nu}),\\
\Pl_{5}(n,\mu) 
      & = U_{\mu}(n-\hat{\mu}) U_{\nu}(n) U^{\dagger}_{\mu}(n+\hat{\nu}-\hat{\mu})U^{\dagger}_{\nu}(n-\hat{\mu}),\\ 
\Pl_{6}(n,\mu)
      & = U_{\mu}(n-\hat{\mu}-\hat{\nu}) U_{\nu}(n-\hat{\nu}) U^{\dagger}_{\mu}(n-\hat{\mu})U^{\dagger}_{\nu}(n-\hat{\mu}-\hat{\nu}),\\
\Pl_{7}(n,\mu)
      & = U_{\mu}(n-2\hat{\nu}) U_{\nu}(n-2\hat{\nu}+\hat{\mu}) U^{\dagger}_{\mu}(n-\hat{\nu})U^{\dagger}_{\nu}(n-2\hat{\nu}),\\
\Pl_{8}(n,\mu)
      & = U_{\mu}(n-\hat{\nu}+\hat{\mu}) U_{\nu}(n-\hat{\nu}+2\hat{\mu}) U^{\dagger}_{\mu}(n+\hat{\mu})U^{\dagger}_{\nu}(n-\hat{\nu}+\hat{\mu}). 
\end{align*}
Then by using the vectors $f(n,\mu)$ defined in 
\eqref{eq:force_fermion}
we obtain the $\C_{FG}$ piece of the force-gradient term given by
% \footnotesize
  \begin{equation*}
  \begin{aligned}
{ \C_{FG}(n,\mu)} & { = 2 \beta
 \left[ \left( 
f(n,\mu) + f(n+\hat{\mu},\nu) - f(n+\hat{\nu},\mu) - f(n,\nu) \right) \operatorname{Re}(\Pl_{1}) \right. } \\
 & { \left. + \left( f(n,\mu) - f(n+\hat{\mu}-\hat{\nu},\nu) - f(n-\hat{\nu},\mu) + f(n-\hat{\nu},\nu) \right) \operatorname{Re}( \Pl_{2}) 
 \right].}
\end{aligned}
\end{equation*}
%\normalsize
{The second derivative of the fermion action is}
\begin{eqnarray}
&& {
\frac{\partial^{2} S_{F}}{\partial q_{\nu}(m) \partial q_{\mu}(n)} = } 
\nonumber \\
&& {
2 \operatorname{Re} \chi^\dagger \left[
\frac{\partial D}{\partial q_\nu(m)} D^{-1} \frac{\partial D}{\partial q_\mu(n)} +
\frac{\partial D}{\partial q_\mu(n)} D^{-1} \frac{\partial D}{\partial q_\nu(m)} -
\frac{\partial^2 D}{\partial q_\nu(m) \partial q_\mu(n)} \right] \xi + }
\nonumber \\
&& {
2 \operatorname{Re} \chi^\dagger \frac{\partial D}{\partial q_\mu(n)}
(D^\dagger D)^{-1} \frac{\partial D^\dagger}{\partial q_\nu(m)} \chi \,, } \nonumber \\
&& {
= 2 \operatorname{Re} \left[ z_{1,m,\nu}^\dagger \frac{\partial D}{\partial q_{\mu}(n)} \xi + 
\chi^\dagger \frac{\partial D}{\partial q_{\mu}(n)} D^{-1} w_{2,m,\nu} - 
\chi^\dagger  \frac{\partial^2 D}{\partial q_{\nu}(m) \partial q_{\mu}(n)} \xi + 
\chi^\dagger \frac{\partial D}{\partial q_{\mu}(n)} D^{-1} z_{1,m,\nu}\right] } \nonumber \\
&& { = 2 \textrm{Re} \left[ z_{1,m,\nu}^\dagger w_{2,n,\mu} + 
w_{1,n,\mu}^\dagger z_{2,m,\nu} - 
\chi^\dagger  \frac{\partial^2 D}{\partial q_{\nu}(m) \partial q_{\mu}(n)} \xi \right] }
\end{eqnarray}
{in terms of the vectors $\chi$ and $\xi$ defined in
\eqref{eq:chi_xi_vectors}.}
Now the fields $z_{1,m,\nu}$ and $z_{2,m,\nu}$ are given by
\begin{equation*}
\begin{aligned}
 z_{1,m,\nu}&:=  D^{\dagger ^{-1}} \frac{\partial D^\dagger}{\partial q_{\nu}(m)} \chi  = D^{\dagger ^{-1}} w_{1,m,\nu} \\
 z_{2,m,\nu}&:=  D^{-1}( \frac{\partial D}{\partial q_{\nu}(m)} \xi + z_{1,m,\nu}) = D^{-1}( w_{2,m,\nu} + z_{1,m,\nu}) 
\end{aligned}
\end{equation*}
with
\begin{equation*}
\begin{aligned}
 w_{1,m,\nu}(n)&: = 
\sum_{n'}\frac{\partial D^\dagger_{n,n'}}{\partial q_{\nu}(m)} \chi(n') = 
\delta_{n,m+\hat{\nu}} \frac{i}{2} (1-\sigma_{\nu})  U^\dagger_{\nu}(m) \chi(m) 
- \delta_{n,m} \frac{i}{2} (1+\sigma_{\nu})  U_{\nu}(m) \chi(m+\hat{\nu}) \,,\\
 w_{2,m,\nu}(n)&: =  
\sum_{n'}\frac{\partial D_{n,n'}}{\partial q_{\nu}(m)} \xi(n') = 
-\delta_{n,m} \frac{i}{2}  (1-\sigma_{\nu}) U_{\nu}(m) \xi(m+\hat{\nu}) 
+ \delta_{n,m+\hat{\nu}} \frac{i}{2} (1+\sigma_{\nu}) U^{\dagger}_{\nu}(m) \xi(m) 
\,.  
\end{aligned}
\end{equation*}

In order to calculate $C_{GF}$ and $C_{FF}$ it is possible to perform the summation of $\sum_{m,\nu}$ before the inversions of $D$ and $D^\dagger$ to get
$z_1$ and $z_2$ which save $\mathcal{O}(V)$ additional inversions for the force gradient terms. It follows for the force gradient term $C_{FF}$
 \begin{equation}\label{eq:CFF}   
 \C_{FF}(n,\mu)=  4 \textrm{Re} \left[ Z_1^\dagger w_{2,n,\mu} + 
w_{1,n,\mu}^\dagger Z_2 - \chi^\dagger  \frac{\partial^2 D}{\partial q_{\mu}(n) \partial q_{\mu}(n)} \xi \cdot f(n,\mu) \right]
 \end{equation}
with 
\begin{equation} \label{eq:Z12}
\begin{aligned}
 Z_1 &:=  D^{\dagger ^{-1}} \sum\limits_{m=1,\nu=1}^{V,2} \left(w_{1,m,\nu} \cdot  f(m,\nu)\right) \,,  \\
 Z_2 &:=  D^{-1} \left(\sum\limits_{m=1,\nu=1}^{V,2} [w_{2,m,\nu} \cdot  f(m,\nu)] + Z_1 \right) \,.
\end{aligned}
\end{equation}
The expression for $C_{GF}$ can be obtained from the one for $C_{FF}$ by
replacing in \eqref{eq:CFF} and \eqref{eq:Z12}
the vector $f$ with $\beta g$ defined in \eqref{eq:force_gauge}.

It is important to mention that the computationally most demanding part of the numerical integration of the Schwinger model 
and quantum field theory in general  is the inverse of the Dirac operator $D^{-1}$. Every momenta update, which includes 
 fermion action \eqref{eq:force_fermion} requires 2 inversions 
of the Dirac operator, the addition of the force-gradient term $\C$  requires 4 more inversions. Therefore leap-frog 
based methods \eqref{eq:leap} and \eqref{eq:leap_multi} need 4 computations of $D^{-1}$ per  time step; schemes \eqref{eq:5stage} 
and \eqref{eq:5stage_multi} 6 times; force-gradient based methods 8 for \eqref{eq:5stage_multi_fg} and \eqref{eq:ad_5stage_multi_fg},
10 for \eqref{eq:5stage_fg} and the 11 stage method \eqref{eq:11_stage} has 12 inversions of the Dirac operator. 
Since we use the multi-rate approach for schemes \eqref{eq:5stage_multi}, \eqref{eq:5stage_multi_fg} and
\eqref{eq:ad_5stage_multi_fg}, which leads generally to fewer macro time steps 
needed than for the standard schemes    
we expect the integrator \eqref{eq:ad_5stage_multi_fg}  will be the 
most efficient choice among the methods considered. 
In the next section we present numerical tests of this prediction.

\section{Numerical results }
In this section we apply the numerical integrators \eqref{eq:leap} -- \eqref{eq:ad_5stage_multi_fg} to compute the molecular dynamics step for 
the Schwinger model \eqref{eq:schwinger} when studied with the HMC algorithm.  
We consider a $32$ by $32$ lattice with a coupling constant $\beta= 1.0$
and mass $m_0=-0.231367$. The parameters were taken from \cite{ChJaNaPo05} and correspond
to the scaling variable $z=0.2$ defined in \cite{ChJaNaPo05}.We have chosen them to simulate close to 
the scaling limit with light fermions and also to
increase the impact of the fermion part of the action. We use one thermalised gauge configuration. For each integrator and value of the step-size we generate
 $200$ independent sets of momenta and integrate the equations of motion on a trajectory of length $\tau=2.0$.
 We compute the absolute error $|\Delta H|$ and estimate its statistical error from the standard deviation. Also the parameter $M$ is chosen in such a way to make 
 micro step size to be $10$ times smaller than the macro step size $h$. 
 
\begin{figure}[h!]
\centering
\includegraphics[width=0.6\linewidth]{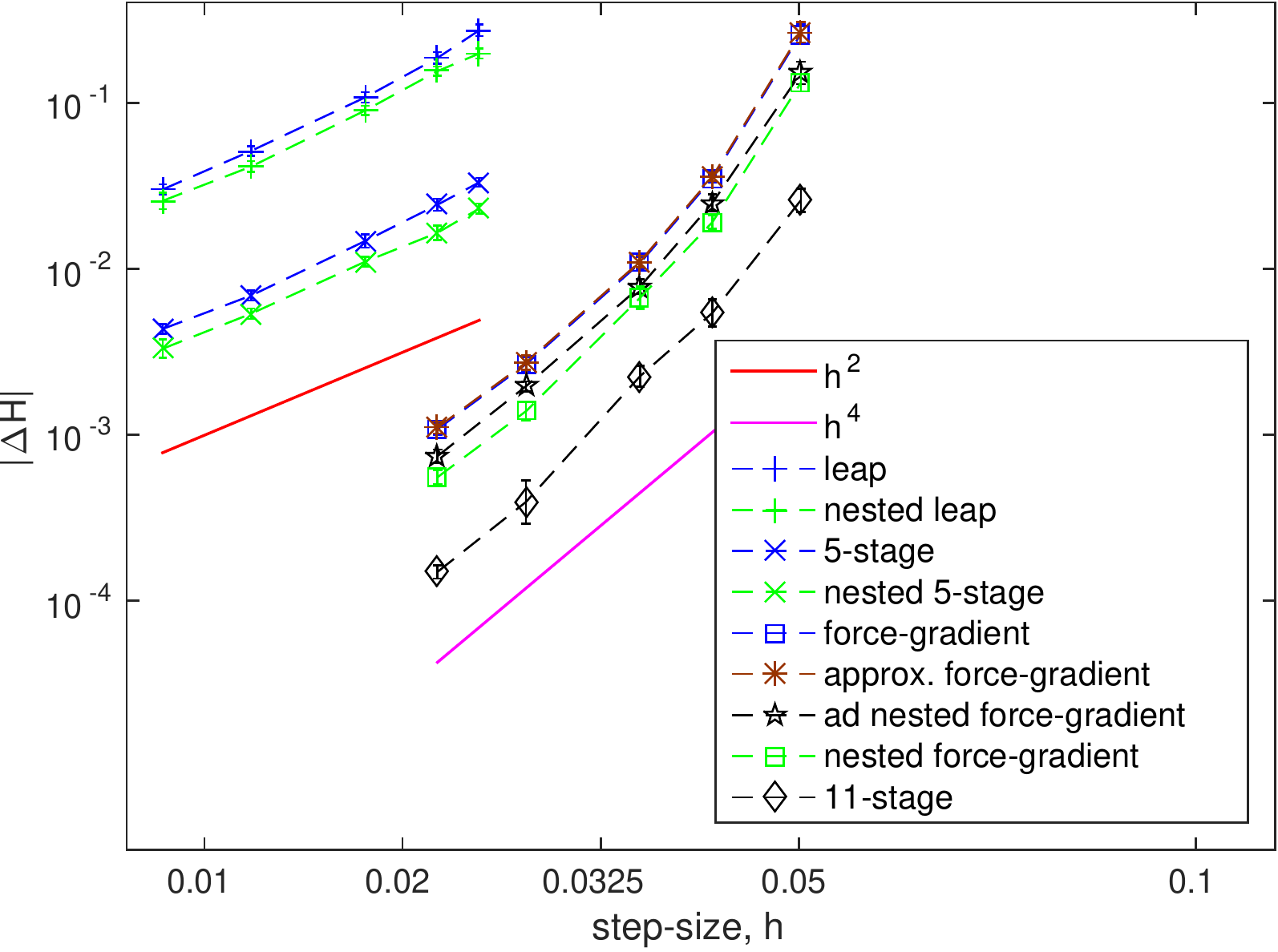}  \vspace{-9pt}
  \caption{Absolute error for different numerical integrators.}
\label{fig:1}
\end{figure}
Figure \ref{fig:1} presents the comparison between the numerical 
integrators \eqref{eq:leap} -- \eqref{eq:ad_5stage_multi_fg}. It  shows the absolute error $|\Delta H|$ versus the step-size 
of the numerical scheme. 
Here the multi-rate schemes  \eqref{eq:leap_multi}, \eqref{eq:5stage_multi}, \eqref{eq:5stage_multi_fg}
and \eqref{eq:ad_5stage_multi_fg} outperform their standard versions as 
expected. Also it is easy 
to see  that the scheme \eqref{eq:11_stage} has the best accuracy and  the nested force-gradient method \eqref{eq:5stage_multi_fg} 
just slightly edges  the adapted nested force-gradient scheme \eqref{eq:ad_5stage_multi_fg}.  

\begin{figure}[h!]
\centering
\includegraphics[width=0.6\linewidth]{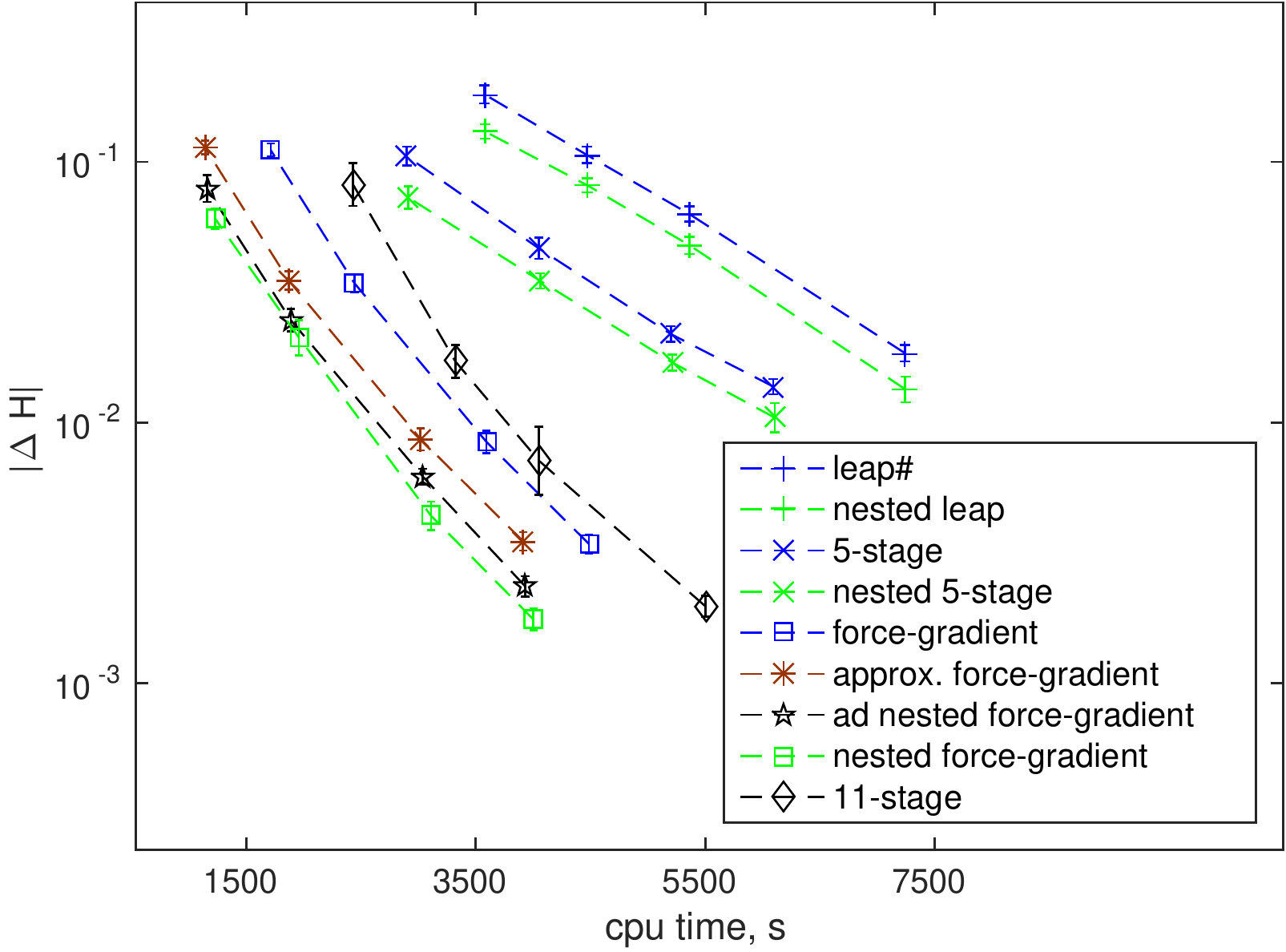}  \vspace{-9pt}
  \caption{Computational costs for different methods}
\label{fig:2}
\end{figure}

Figure \ref{fig:2}  presents the CPU time, required for the proposed integrators \eqref{eq:5stage}--\eqref{eq:ad_5stage_multi_fg}, 
versus the achieved accuracy.  We can observe that the nested force-gradient method \eqref{eq:5stage_multi_fg} and adapted nested force-
gradient method \eqref{eq:ad_5stage_multi_fg} show much better results in terms of a computational efficiency 
than  the integrators \eqref{eq:5stage_multi} and  \eqref{eq:5stage_fg}; and even compared to the 11 stage scheme 
\eqref{eq:11_stage}. Here we can see that the modification of \eqref{eq:5stage_fg} proposed in \cite{YiMa} also performs
better than its original version. It shows almost similar computational costs as  nested versions of the force-gradient approach \eqref{eq:5stage_multi_fg} -\eqref{eq:ad_5stage_multi_fg}, 
since it has the same number of $D^{-1}$~(see Table \ref{tab:1}). But it is less efficient because the proposed nested approach is more precise.
%The numerical scheme \eqref{eq:ad_5stage_multi_fg} slightly faster then the scheme \eqref{eq:5stage_multi_fg}.
\vspace{-9pt}
\begin{center}
\begin{table}[h!]
\begin{tabular}{||c|c|c|c|c||}
\hline
 Integrator:                        & step size $h$ & $M$  & $D^{-1}$ per step & $D^{-1}$ per trajectory\\ \hline
 5 stage method                     &    0.0294     & -    & 6                 & 420\\
 nested 5 stage method              &    0.0286     &  700 & 6                 & 408 \\
 5 stage force-gradient             &    0.0550     & -    & 10                & 370 \\
  approx.~force-gradient \cite{YiMa}  &    0.0540     & -    & 8                 & 290  \\
 nested  force-gradient             &    0.0560     &  450 & 8                 & 285  \\
adapted  nested force-gradient      &    0.0560     &  450 & 8                 & 285  \\
 11 stage method                    &    0.0625     & -    & 12                & 384\\
 \hline
\end{tabular} 
\caption{Step-sizes and number of inversions of $D$ per step and per trajectory  for acceptance rate of 90\%}
\label{tab:1}
\end{table}
\end{center}

Table  \ref{tab:1} shows the number of  inversions of the Dirac operator $D$, which is needed to reach  90\% acceptance rate of the HMC. 
Since $D^{-1}$ is the most computationally  demanding part it is important to see how many of these inversions are required per 
each trajectory.  From Table  \ref{tab:1} it easy to see that the adapted nested force-gradient method
\eqref{eq:ad_5stage_multi_fg}  and nested 
force-gradient method  \eqref{eq:5stage_multi_fg}  
need the  least number of $D^{-1}$  per trajectory to reach the chosen acceptance rate $ \approx 90 \%$. 

We can also claim that methods \eqref{eq:5stage_multi_fg} and \eqref{eq:ad_5stage_multi_fg}   have a potential to perform 
even better with respect to the 
computational effort in the case of lattice QCD problems, since the impact of the fermion action \eqref{eq:fer_action} 
and the computational time  to obtain  the inversion of the Dirac operator $D$ is much more significant.

%\begin{figure}[h!]
%\centering
%%\includegraphics[width=0.8\linewidth]{p_vs_d_inv}  \vspace{-10pt}
%  \caption{Number of $D^{-1}$ per trajectory}
%\label{fig:3}
%\end{figure}

%Figure \ref{fig:3} shows how many inversions of the Dirac operator $D$  are needed to reach  a certain accuracy rate of the HMC. 
%Since $D^{-1}$ is the most computationally  demanding part it is important to see how many of these inversions are required per trajectory.  From 
%Figure \ref{fig:3} it is easy to see the adapted nested force-gradient method
%\eqref{eq:ad_5stage_multi_fg}  and nested 
%force-gradient method  \eqref{eq:5stage_multi_fg}  
%need the  least number of $D^{-1}$  per trajectory to reach the chosen accuracy of $ \approx 90 \%$.  

%%%%%%%%%%%%%%%%%%%%%%%%%%%%%%%%%%%%%%%%%%%%%%%%%%%%%%%%%%%%%%%%%%%%%%%%%
\section{Conclusions and outlook}

We presented the nested  force-gradient approach \eqref{eq:5stage_multi_fg} and its adapted version \eqref{eq:ad_5stage_multi_fg} applied to a model problem in 
quantum field theory, the two-dimensional Schwinger model. The derivation of 
the force-gradient terms was given and the Schwinger model was introduced. 
Nested force-gradient schemes
%\eqref{eq:ad_5stage_multi_fg} 
seem to be an optimal choice with relatively 
high convergence order and low computational effort. Also it would be possible to improve the algorithm \eqref{eq:ad_5stage_multi_fg} by measuring
the  Poisson brackets of the shadow Hamiltonian of the proposed integrator
and then tuning the set of optimal parameters, e.~g. micro and  macro step sizes.\\

In  future work we will apply this approach to the HMC algorithm for numerical 
integration in Lattice QCD.  
Here we expect the adapted nested-force gradient scheme to outperform the original one, if we further partition 
the action into more than two parts, by using techniques to factorize the fermion determinant: 
less force-gradient information is needed for the most expensive action, and only leap-frog steps are needed for 
the high frequency parts of the action.

%Due to the  
%universal formulation of the proposed scheme, the most challenging part will be  to derive the force-gradient term and  
%to implement 
%it in the available source code. 

%%%%%%%%%%%%%%%%%%%%%%%%%%%%%%%%%%%%%%%%%%%%%%%%%%%%%%%%%%%%%%%%%%%%%%%%%

\section*{Acknowledgments}
This work is part of project B5 within the SFB/Transregio 55
{\em Hadronenphysik mit Gitter-QCD} funded by DFG (Deutsche Forschungsgemeinschaft).

%%%%%%%%%%%%%%%%%%%%%%%%%%%%%%%%%%%%%%%%%%%%%%%%%%%%%%%%%%%%%%%%%%%%%%%%%

\end{document}